\documentclass[11pt,a4paper]{article}
\usepackage[left=3cm, right=3cm, top=3cm, bottom=3cm]{geometry}

\usepackage{graphicx}
\usepackage{subfig}
\usepackage{mathtools}
\usepackage{float}
\usepackage{amssymb}
\usepackage{amsthm}
\usepackage{thmtools}
\usepackage{xcolor}
\usepackage{nameref}
\usepackage{hyperref}

\usepackage{authblk}

\DeclareMathOperator{\diff}{d}
\newcommand{\dt}{\mbox{d}}
\newcommand{\N}{\mathcal{N}}

\title{\large \bf   Data assimilation for the  stochastic Camassa-Holm equation using  particle filtering: a numerical investigation}
\date{}

\author{Colin J. Cotter}
\author{Dan Crisan}
\author{Maneesh Kumar Singh\thanks{Corresponding author: maneesh-kumar.singh@imperial.ac.uk}}
\affil{Department of Mathematics,
Imperial College London, UK}

\begin{document}
	\maketitle

\begin{abstract}
    In this study, we explore data assimilation for the Stochastic Camassa-Holm equation through the application of the particle filtering framework. Specifically, our approach integrates adaptive tempering, jittering, and nudging techniques to construct an advanced particle filtering system. All filtering processes are executed utilizing ensemble parallelism. We conduct extensive numerical experiments across various scenarios of the Stochastic Camassa-Holm model with transport noise and viscosity to examine the impact of different filtering procedures on the performance of the data assimilation process. Our analysis focuses on how observational data and the data assimilation step influence the accuracy and uncertainty of the obtained results.
\end{abstract}

\section{Introduction}\label{sec1:intro}

Data assimilation (DA) is a set of methodologies that integrate past knowledge, represented as
numerical models of a system, with newly acquired observational data from the same
system. This tool is used in many domains, including meteorology, oceanography, and
environmental research, to merge observational data with numerical models to improve
forecast and simulation accuracy. A concise overview of DA is presented in \cite{asch2016data}, including the references contained within. For stochastic systems, data assimilation can be rigorously formulated as stochastic filtering. In this work, we emphasise a stochastic filtering problem where a hidden stochastic process (signal) is observed at discrete times with noise. The nonlinear filtering problem consists of computing the law of the signal, given observations that are collected sequentially. More details on stochastic filtering can be found in \cite{bain2009fundamentals} and references therein. 

In this study, we investigate the data assimilation method for a nonlinear stochastic partial differential equation that corresponds to a viscous shallow water equation. In particular, we examine the particle filter methodology for a stochastic Camassa–Holm (CH) model with transport noise (referred to as Stochastic Advection by Lie Transport, or SALT).  The deterministic CH equation \cite{camassa1993integrable} admits solutions with singularities, which possess a sharp peak at the apex of their velocity profile. By following the variation principle approach in stochastic fluid dynamics \cite{holm2015variational}, the stochastic Camassa–Holm (SCH) equation in the SALT framework is derived in \cite{crisan2018wave}. The formation of  Peakon (peaked soliton solutions) of the SCH model due to stochastic transport is investigated in \cite{bendall2021perspectives}. This model is useful in providing a 1+1 SPDE example where the behaviour of SALT in data assimilation can be easily investigated. However, there is a disadvantage: the solutions become rougher as time progresses. Numerical solutions behave like a space-time random field at longer times, leading to an atypical data assimilation problem (which is easy to solve numerically since the long-time solution of the filtering problem is a steady state distribution describing this random field). In this paper we incorporate a viscous dissipation to the equations to control the regularity of the solution at longer times; this provides a more challenging and informative data assimilation problem to benchmark our filtering methods against.

The development of ensemble-based data assimilation, an alternative to variational data assimilation, has received greater research interest in recent years. A detailed discussion on both approaches can be found in \cite{reich2015probabilistic,van2019particle}. In \cite{llopis2018particle}, particle filtering is discussed for the stochastic Navier–Stokes equation with linear additive type noise. The data assimilation method using a particle filter for the two-dimensional Euler equation and quasi-geostrophic model is investigated in 
\cite{cotter2020particle,cotter2020data}. A tempering-based adaptive particle ﬁlter to infer from a
partially observed stochastic rotating shallow water (SRSW) model is studied in  \cite{lang2022bayesian}.
Recently, a  lagged particle filter has been introduced for stable filtering of high-dimensional state-space models in \cite{ruzayqat2022lagged}.

The organization of this paper is as follows. In section \ref{sec2}, we introduce the stochastic  Camassa–Holm equation and discuss the numerical approximation of the model problem. In section \ref{sec3}, we briefly review the bootstrap particle filter and its consequent augmentations with jittering, tempering and nudging procedures. In section \ref{sec4}, a numerical study is conducted, illustrating the behaviour of data assimilation methods. Finally, we summarize with some concluding remarks in section \ref{sec5}.

\section{The SCH model and its discretisation}\label{sec2}

The deterministic viscous  Camassa–Holm is expressed as the following evolution equation,
\begin{equation}\label{CHe1}
 \begin{array}{ll}
 m = u - \alpha^2 \partial_{xx}u, \\[4pt]
\dt m -\mu \partial_{xx}m  \, \dt t+ (\partial_{x}m + m \partial_{x})\diff v = 0, \quad  \mbox{with} \,\, \diff v = u \dt t \\[4pt]
u(x,0) = u^{0}(x),\quad m(x,0) = u^{0}(x) - \partial_{xx}u^{0}(x), \,\, x \in [0,L],
\end{array}   
\end{equation}
for the evolution in time $t \in [0,T)$
of the fluid momentum density $m(x,t)$ and the velocity $u(x,t)$
solved on the spatial domain  $[0,L]$ with periodic boundary conditions. 

In the SALT framework, this deterministic PDE is transformed into a stochastic PDE by changing $\diff v$ to
\begin{equation}
\diff v = u \dt t + \diff U,
\end{equation}
where $U$ is some stochastic process.
Instead of the noise expansions of the form $\diff U = \sum_{k=1}^{K} \xi^{k} \dt W^{k}$
used in some previous works on SALT, in this paper, we use Gaussian random space-time fields obtained from the
Mat{\'e}rn formula,
\begin{equation}\label{noise_solve}
    (I - \kappa^{-2} \nabla^{2})^{k}\dt U(x,t) = \eta \dt W(x,t),
\end{equation}
where $\diff W(x,t)$ is a space-time white noise that 
is cylindrical in space and Ito in time.
The coefficients $\eta$ and $\kappa$ determine the expected smoothness of the process $U$. For this work, we consider $\eta=1$, and $k=3$.

Henceforth, we consider the SCH equation,
\begin{equation}\label{schce1}
\begin{array}{ll}
 m = u - \alpha^2 \partial_{xx}u, \quad  t \in (0,T], \,\, x \in (0,L)\\[4pt]
\dt m -\mu \partial_{xx}m \, \dt t+ (\partial_{x}m + m \partial_{x})\diff v = 0, \quad  \mbox{with} \,\, \diff v = u \dt t+ \dt U \\[4pt]
u(x,0) = u^{0}(x),\quad m(x,0) = u^{0}(x) - \partial_{xx}u^{0}(x), \,\, x \in [0,L].
\end{array}
\end{equation}

For spatial discretisation, we will be interested in approximating solutions of the SCH that are periodic functions in the spatial variable. We will consider a finite element discretisation on a uniform mesh of
the interval $I=[0,L]$ with $N$ cells of width $h=L/N$.

The space-time white noise $\diff W(x, t)$ is approximated by $\diff W_h(x, t)$, defined
as 
\begin{equation}
\diff W_h(x, t) = \sum_{i=1}^N \frac{1}{A_i^{1/2}}\phi_i(x)\diff W_i(t),
\end{equation}
where $A_i$ is the width of cell $i$ (which is equal to $h$ for a uniform grid), $\phi_i$ is the indicator 
function of cell $i$, and $\{W_i(t)\}_{i=1}^N$ are $N$ iid Brownian motions. In other words, $\int_{t_A}^{t_B}\diff W_h \in Q_h$ for any times $t_A<t_B$, where $Q_h$ is the piecewise constant
finite element space (DG0).
This is a low-order spatial approximation of the space-time white noise; higher-order approximations can be obtained but the required square root factorisation of the finite element mass matrix is not diagonal and we did not implement this. 
\cite{croci2018efficient} provided an efficient formulation for continuous finite element fields using 
a generalised form of the square root that exploits the local assembly procedure.

We then approximate $\diff U$ with $\diff U_h:= \diff U_k$, where $\{U_j\}_{j=1}^k$ are space-time Gaussian random fields in the continuous linear Lagrange (P1) 
finite element space on the mesh, denoted here as $V_h$.
Then we approximate the solution of \eqref{noise_solve}
by solving $k$ approximated second-order elliptic problems, according to
\begin{equation}
\label{noise_solve2}
      (\Delta U_j, v) + \kappa^{-2}(\nabla \Delta U_j, \nabla v) =
      \left\{
      \begin{array}{r l}
      \eta(\Delta W_h, v), & j=1 \\
      (\Delta U_{j-1}, v), & j>1 \\
      \end{array}\right.,
      \quad \forall v \in V_h, \quad
      j=1,2,\ldots, k,
\end{equation}
where $(\cdot,\cdot)$ represents the usual $L^2$
inner product for functions on $I$,
$\Delta U_j=\int_{t_A}^{t_B} \diff U_j(t)$,
$\Delta W_h = \int_{t_A}^{t_B} \diff W_h(t)$, with 
the formula holding for any $t_A<t_B$, and $\Delta U_j,\Delta W_h$, for $j=1,2,\ldots, k$.

Eq. (\ref{noise_solve2}) is solved numerically and the solution $\dt U_{3}$ is added in the velocity term $ v = u \dt t+ \dt U$ to obtain our stochastic Camassa-Holm (SCH) model.

Then, the semidiscrete numerical scheme seeks
$m(t)\in V_h$ and $u(t)\in V_h$ such that
\begin{equation}\label{sde1}
\begin{array}{ll}
	(u,\psi) + \alpha^2 ( \partial_{x}u,\partial_{x}\psi) - (m,\psi) =0, \quad \forall \psi \in V_{h},\\[4pt]
	 (m_t, \phi) + \mu ( \partial_{x}m, \partial_{x}\phi) + (m\partial_{x}v, \phi) - (mv, \partial_{x}\phi)=0, \quad \forall \phi \in V_{h}.
\end{array}
\end{equation}
This is equivalent to a finite dimensional stochastic
differential equation on $\mathbb{R}^N$.
More details on the selection of the initial conditions are given in a later section.

For the time discretization, we select a uniform time step $\Delta t = T/M$ and $t_{n} = n \Delta t, \,\,n=1,2,\ldots, N_{T}$, and solve for 
$w^n\approx w(x,t_n)$ and $u^n\approx u(x, t_n)$
Anticipating an explicit time discretisation of the noise term, we define $\Delta U^n_h$ by taking $t_A=t_n$, $t_B=t_{n+1}$ in the definition of $\Delta U$ above.

Then we use $\Delta U^n_h$ in a midpoint rule discretisation (leading to a Stratonovich method), 
according to finding $w^{n+1},u^n \in V_h$
\begin{equation}\label{fde1}
\begin{array}{ll}
	(u^{n+1},\psi) + \alpha^2 ( \partial_{x}u^{n+1},\partial_{x}\psi) - (m^{n+1},\psi) =0, \quad \forall \psi \in V_{h},\\[4pt]
	 (m^{n+1}-m^{n}, \phi) + \mu \Delta t ( \partial_{x} m^{n+1/2}, \partial_{x}\phi) \\[4pt]
 \qquad  + (m^{n+1/2}\partial_{x}v^{n+1/2}, \phi) - (m^{n+1/2} v^{n+1/2}, \partial_{x}\phi)=0, \quad \forall \phi \in V_{h},
\end{array}
\end{equation}
where $\Delta v^{n+1/2}=\Delta t u^{n+1/2}+\Delta U^n_h$ 
and $m^{n+1/2}=(m^{n+1}+m^n)/2$, etc.

\section{Data assimilation methods}\label{sec3}

In section \ref{sec2}, we defined the SPDE providing the unknown signal, i.e., the system we are interested in performing Bayesian inference upon.
In this work, we will use the language
of stochastic filtering to provide the background framework for a
Bayesian inference case study for the SCH model. 

Let $X$ and $Y$ be two processes defined on the probability space $(\Omega, \mathcal{F}, \mathbb{P})$. The process $X$ is usually called the signal process or the truth and $Y$ is the observation process. The pair of processes $(X, Y)$ forms the basis of the nonlinear filtering problem: find the best approximation of the posterior distribution of the signal $X_t$, denoted by $\pi_t$ given the observations $Y_{1}, Y_{2},\ldots, Y_t$. In our context, the observations
consist of noisy measurements of the true state recorded at discrete times and they are taken at locations on a data grid $\mathcal{G}_{d}$, defined later. The data assimilation is performed at these times, which we call the assimilation times.

In this work, we discuss the approximation of the posterior
distribution of the signal by \emph{particle filters}. These
sequential Monte Carlo methods generate approximations of the posterior distribution using sets of \emph{particles}, which are empirical samples from the conditional distribution of $X$.
Particle filters are employed to make inferences about the signal process. This involves utilizing Bayes' theorem, considering the time-evolution induced by the signal $X_t$, and taking into account the observation process $Y_t$. The observation data $Y_t$ is, in our case, an $M$-dimensional process  that consists of noisy
measurements of the velocity field $u$ taken at a point belonging to the data grid $\mathcal{G}_d = \{iL/M\}_{i=0}^{M}$:
\begin{equation}\label{obseq}
    Y_t := \mathcal{P}_{d}^{s}(X_t) + V_t,
\end{equation}
where $\mathcal{P}_{d}^{s}: \mathcal{G}_{s} \rightarrow \mathcal{G}_{d}$ is a projection operator from the solution space to the data grid $\mathcal{G}_{d}$ and $V_{t} = \N(\mathbf{0}, I_{\sigma})$, where $I_\sigma = \mbox{diag}(\sigma^{2}_{1}, \sigma^{2}_{2}, \ldots, \sigma^{2}_{M} )$. It is important to note here that $ M 
\neq N$.
While we assumed standard normal distributions for $V_{t}$, the methodology presented is valid for a broader range of distributions. The ensemble of particles evolved
between assimilation times according to the law of the signal.

Next, we explain briefly the various types of particle filters used in this article. Before going into the detail of particle filters, we introduce some technical terms, which are used to explain filters.  The  likelihood weight function is defined as 
\begin{equation}\label{loglikeq}
    \mathcal{W}(\mathbf{X},\mathbf{Y}) = \exp\left(-\dfrac{1}{2}\sum_{i=0}^{M}\left\|\dfrac{\mathcal{P}_{d}^{s}(X_{i})-Y_{i}}{\sigma_{i}}\right\|^{2}_{2}\right)
\end{equation}
with $M$ as the number of observation grid points. To measure the variability
of the weights (\ref{loglikeq}) of particles, we use the \emph{effective sample size},
\begin{equation}\label{ess}
    \mathrm{ESS}(\overline{\mathbf{w}}) = \left(\sum_{n=1}^{N_p}(\overline{\mathbf{w}}_{n})^2 \right)^{-1}, \quad \overline{\mathbf{w}} := \mathbf{w}\left(\sum_{n=1}^{N_p}w_{n} \right)^{-1}, \quad w_{n} = \mathcal{W}(\mathbf{X}^{(n)},\mathbf{Y}^{(n)}),
\end{equation}
where $n$ is the particle index, and $N_p$ is the number of particles.
$\mathrm{ESS}$ quantifies the variance of weights. The $\mathrm{ESS}$ value approaches $N_p$ if the particle weights are almost uniform, and it is close to one when fewer particles have large weights and the remaining particles have small weights. We resample whenever $\mathrm{ESS}$ drops below a specified $N_{p}^{*}$ threshold. For all the numerical experiments, we let $N_{p}^{*} = 80$ to be a threshold.

\subsection{Particle Filter: Basics terminology}
In this subsection, we briefly discuss the basics of particle filters, so that we can present our results in context. We mostly describe the methodology; more details on why it works can be found in the references \cite{llopis2018particle,cotter2020particle,cotter2020data}. 

\subsubsection{Bootstrap Particle Filter}
The bootstrap particle filter is the basic particle filter, also called a Sampling Importance Resampling filter. In the bootstrap filter, given an initial distribution of particles (obtained as samples from a prior distribution for the initial state), each particle is propagated forward according
to the signal equation (the spatially discretised SCH equation in our case), with independent realisations of the Brownian motions. Here, and in the more sophisticated particle filter formulations later, 
we consider intermittent data assimilation in intervals of length $\Delta \tau$, subdivided into model timesteps $\Delta t$
with $\Delta \tau/\Delta t=N_s$ some positive integer. When we discuss discrete time we will use the suffix $u^n, m^n$ to
indicate the solution $n$ timesteps after the last assimilation time, i.e. the index $n$ resets to 0 after the most recent
data has been assimilated, for the purposes of presentation here.

The empirical distribution,
\begin{equation}
\diff \mu_F = \sum_{i=1}^{N_p}\frac{1}{N_p}\delta(X - X_i),
\end{equation}
where $\delta$ is the Dirac measure,
is an approximation of the prior (forecast) distribution for the signal, before receiving the observations.

Subsequently, utilizing partial observations, weights for new particles are calculated.
This is done by computing the likelihood weight function
\eqref{loglikeq} for each particle, and then renormalising so that the complete set of weights sums to 1. The weighted empirical distribution,
\begin{equation}
\diff \mu_A = \sum_{i=1}^{N_p}w_i\delta(X - X_i),
\end{equation}
where $\{w_i\}_{i=1}^{N_p}$ are the normalised weights,
is an approximation of the posterior (analysis) distribution for the signal, conditional on the received observations.

Next, a selection process is used on the weighted particles. This is a statistical procedure that aims to 
find a new equally weighted set of particles that approximate the same distribution as the old nonequally weighted set of particles is approximating.
On average, the particles with larger weights will be duplicated, while the particles which have smaller weights will be eliminated. In the simplest case, this is done by sampling with replacement from the ensemble of particles using the multinomial distribution described by the weights. In this work, we use the systematic resampling algorithm \cite{gordon1993novel},
to reduce the sampling error.

The $\mathrm{ESS}$ is a crude diagnostic that 
measures how far the weights are from being uniform.
The value of $\mathrm{ESS}$ typically drops fast for higher dimensional problems because the sample degenerates rapidly. This is due to the trajectory of the particles quickly diverging from the trajectory of the signal.
As a result, the particles fail to give a better approximation of the posterior distribution. To overcome this situation, one would require a huge number of particles.

To resolve the filter degeneracy,  we replace the direct resampling from the weighted predictive approximations by a \emph{tempering} procedure combined with \emph{jittering} and \emph{nudging} procedures, described below.
  
\subsubsection{Tempering and jittering} 
As discussed above, due to the sample degeneracy, 
the $\mathrm{ESS}$ value will rapidly fall below the tolerance value $N_{p}^{*}$. The purpose of \emph{tempering} is to take incremental steps between
the approximations of the prior and posterior distributions, resampling on each step, to maintain
a high $\mathrm{ESS}$ value. In each tempering 
step $k=1,\ldots,N_\theta$, the particle likelihood weights are evaluated
and scaled by $0<\Delta \theta_k<1$, with $\sum_{k=1}^{N_k}\Delta \theta_k=1$, where
$\Delta \theta_k$ is chosen so that $\mathrm{ESS}>N_p^*$
for that step. In our work we use an adaptive tempering procedure as discussed in \cite{cotter2020particle}. This procedure repeatedly reduces $\Delta \theta_k$ until the condition holds, keeping the number of tempering steps to a minimum.

After each tempering step, the particles are resampled
according to the scaled (and normalised) weights. 
This alone is insufficient to prevent the accumulation of a large number of duplicates in the particle ensemble.
To remove the duplication in the resulting ensembles, we employ \emph{jittering}. This can be understood by using the fact that distributions on the state at assimilation time $\tau_n$ are equivalent to distributions on the joint distribution of the state at assimilation time $\tau_{n-1}$ together with the Brownian increments
from $\tau_{n-1}$ to $\tau_n$. The equivalence comes
because given a sample from the latter, we can solve
the signal equation forwards with the initial condition given by the state value $X_{\tau_{n-1}}$ at time $\tau_{n-1}$, using
the realisation of the Brownian increments $dW$, to obtain
a state value $X_{\tau_n}$ at time $\tau_n$. Each particle
can thus be represented by $(X_{\tau_{n-1}},dW)$, which
may be duplicated after resampling.
Using the disintegration formula $\pi(X_{\tau_{n-1}}, dW)=\pi(X_{\tau_{n-1}})
\pi(dW|X_{\tau_{n-1}})$, we see that keeping $X_{\tau_{n-1}}$ the
same, but obtaining a new $dW$ sample from $\pi(dW|X_{\tau_{n-1}})$, produces another consistent sample from $\pi(X_{\tau_{n-1}},\diff W)$. This is the jittering technique that allows duplicates to be replaced by different samples 
from the same distribution. This is achieved by performing tempering using the $(X_{\tau_{n-1}},dW)$ representation, and only updating to $X_{\tau_{n-1}}$ once
the tempering step is complete. At each tempering 
step, after resampling, the noise realisation $dW$ for
each particle is moved using a Monte Carlo Markov Chain (MCMC) method on the tempered posterior
distribution with the likelihood function scaled
by $\theta_k=\sum_k\Delta \theta_k$ at step $k$. Since
the samples are from the conditional distribution
$\pi(dW|X_{\tau_{n-1}})$, this can be done independently for each particle.

The MCMC method describes a sequence of samples
$dW$ from $\pi(dW|X_{\tau_{n-1}})$, namely $dW^0,dW^1,dW^2,\ldots$, where $dW^0$ is the sample
of $dW$ given after resampling.
Note that here $dW$ represents \emph{all} of the 
Brownian increments from $\tau_n$ to $\tau_{n+1}$. After space and time
discretisation, this 
is a finite array of numbers (of dimension $N_s\times N$ for
our discretisation choices)
whose prior distribution is iid 
$\mathcal{N}(0,\Delta t)$. We use the Preconditioned Crank Nicholson (PCN) algorithm \cite{cotter2013mcmc} to move the particles. PCN proposes
a new sample $\diff W \mapsto (2-\delta)/(2+\delta)\diff W + (8\delta)^{1/2}/(2+\delta)\diff\hat{W}$, where $\diff\hat{W}$
is a new sample from the prior distribution, and $\delta>0$ is a time-stepping parameter. The proposal is accepted with probability
\begin{equation}
a = \max(1, \exp(w_0 - w_1)),
\end{equation}
where $w_0$ is the $\theta_k$-weighted negative log-likelihood 
for $X_n$ obtained from the old $\diff W$ and 
$w_1$ is the same quantity obtained from the new
proposal, otherwise the old value of $\diff W$ 
is repeated. Larger $\delta$ means that the proposal is moved further away and smaller $\delta$ means that the proposal is more likely to be accepted. The number
of MCMC iterations (we call them ``jittering steps")
per tempering steps is fixed. It is not necessary for 
the MCMC algorithm to converge in statistics, just that the duplicated particles are sufficiently spread. The accept-reject
criteria insures statistical consistency as proved in \cite{beskos2014stability}.

\subsubsection{Nudging} 
In the nudging particle filter framework, we 
introduce a time dependent control variable that 
we may choose, so that
particles are ``nudged'' towards regions where the observations suggest the signal is likely to be. This \emph{nudging} is done in a way that preserves the 
consistency of the particle filter, \emph{i.e.} the
particle ensemble remains a consistent set of samples
from the prior distribution, after appropriate 
modification of the weights. Presently, we justify
this at the level of the spatial semidiscretisation,
which is interpretable as an SDE, of the form
\begin{equation}
\diff x = f(x)\diff t + G(x)\diff W, \, x(0)=x_0,
\end{equation}
where $x\in \mathbb{R}^N$, $f:\mathbb{R}^N\to \mathbb{R}^N$, $G: \mathbb{R}^N\to \mathbb{R}^{N\times Q}$ (\emph{i.e.}, $G(x)$ is an $N\times Q$ matrix for 
each $x\in \mathbb{R}^N$), and $W(t)$ is a $Q$-dimensional Brownian motion. If one instead solves
the modified SDE,
\begin{equation}
\diff \hat{x} = f(\hat{x})\diff t + G(\hat{x})(\lambda(t)\diff t + \diff W),\, \hat{x}(0) = x_0,
\end{equation}
where $\lambda(t) \in \mathbb{R}^Q$, then the joint 
probability measure for $\hat{x}(t)$ for $t\in (0,T)$ is
absolutely continuous with respect to the joint 
probability measure for $x(t)$ on the same range, and
the Radon-Nikodym derivative from one measure to the
other is
\begin{equation}
\label{eq:girsanov}
G = \exp\left(
-\int_0^T \frac{1}{2}|\lambda(t)|^2 \diff t 
+ \int_0^T \lambda(t)\cdot\diff W
\right).
\end{equation}
This means that we can choose $\lambda$ to reduce
the likelihood weight, but we must pay the price of 
multiplying this weight by $G$. 

In particular, in this work we use the Girsanov formula to correct the solution of the spatially discrete SCH equation \eqref{sde1} to keep the particles closer to the true state. To implement this, we update the SCH model (\ref{schce1}) with a ‘nudging term’, replacing the $j=1$ case of \eqref{noise_solve} with
\begin{equation}
\label{eq: SCH nudge}
      (\Delta U_j, v) + \kappa^{-2}(\nabla \Delta U_j, \nabla v) =
      \eta(\Delta W_h + \Delta \Lambda, v),\quad\forall v\in V_h,
\end{equation}
where $\Delta\Lambda=\Lambda(t)(t_B-t_A)$, and 
$\Lambda(t)\in Q_h$. The Girsanov formula \eqref{eq:girsanov} can be rewritten as 
\begin{equation}
\label{eq:girsanov DG}
G = \exp\left(
-\int_0^T \frac{1}{2} (\Lambda^2(t), 1/A) \diff t 
+ \int_0^T (\Lambda(t), \diff W_h/A)
\right),
\end{equation}
where $A\in Q_h$, such that $\int_{e_i}A \diff x=1$,
for each cell $e_i$. In our discrete time approximation
of this, we use $\Delta U_h=\Delta U_k$ from
\eqref{eq: SCH nudge} with $t_A=t_n$ and $t_B=t_{n+1}$
in \eqref{fde1}, and we use the approximated time integral,
\begin{equation}
\label{eq:girsanov DG approx}
G \approx G_\Delta t = \exp\left(
-\sum_{n=0}^{N_s}\left(\frac{1}{2} ((\Lambda^n)^2, 1/A) \Delta t 
+ (\Lambda^n, \Delta W_h^n/A)
\right)\right).
\end{equation}

After choosing $\lambda(t)$, the particles will have a new
weights according to Girsanov’s theorem, given by
\[
\tilde{W}(u,Y,\Lambda) = W(u,Y) + G(\Lambda).
\]
Hence, it makes sense to maximise $\tilde{W}$ over $\Lambda$.
However, a critical aspect is that Girsanov's theorem only
holds if $\lambda(t)$ only depends on $W(s)$ for $s<t$.
Our strategy is to incrementally optimise $\lambda(t)$
as we reveal $W(t)$; we have to adapt $\lambda(t)$ to past noise. After time discretisation, this means that we first initialise $\Lambda^n=0$ and $\Delta W^n=0$ for $n=0,1,\ldots,N_s$. Then we optimise $\tilde{W}$ over $\Lambda^0$, keeping the other values of
$\Lambda$ fixed. Next, we randomly sample $\Delta W^0_h$
from the $Q$-dimensional distribution $\N(0,\Delta t I)$.
Then we optimise $\tilde{W}$ over $\Lambda^1$, keeping the other values fixed, and randomly sample $\Delta W^1_h$ from $\N(0, \Delta t I)$, and so on, until we reach
$n=N_s$. This approach was investigated in \cite{cotter2020data} in application to data assimilation for a quasigeostrophic ocean model, but $\Lambda$ was only nonzero in the final stage of the splitting method in the last timestep $n=N_s$, which allows the optimisation problem to be solved using linear least squares. This simplification was due to 
the lack of an available nonlinear optimisation algorithm within that code framework. This is something that we have addressed in the current work, as described below.

In general these optimisation problems are nonlinear, since the observations at $\tau_n$ depend on the entire history of 
$\diff W$ from $t=\tau_{n-1}$ to $\tau_n$ through the nonlinear SDE. We solve these problems numerically, using an gradient descent algorithm (BFGS), with the gradient of the functional computed using the adjoint technique. In our implementation, this is automated using Firedrake \cite{FiredrakeUserManual}, which is built according to the methodology of \cite{farrell2013automated}.\footnote{Since we need to solve repeated optimisation problems with different data, this required minor extensions to Firedrake, namely the {\ttfamily derivative\_components} argument
to {\ttfamily adjoint.ReducedFunctional}, which allows to ignore
derivatives with respect to the observed data in the
minimisation calculation; these changes are now in the main branch of Firedrake.}

After applying the nudging step, it may still be necessary to use tempering and jittering as above.
In that case, the $\theta_k$ adjusted weight formulas described above need to be modified by replacing
$W$ with $\tilde{W}$, as described in \cite{cotter2020data}.
	
\subsubsection{Ensemble parallelism}

\begin{figure}
	\centering
	\includegraphics[scale=0.45]{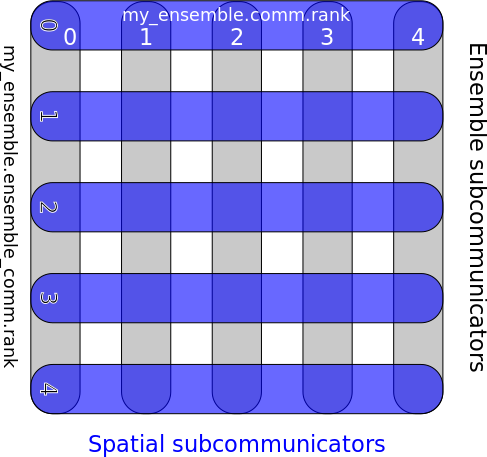}
	\caption{Spatial and ensemble parallelism for an ensemble with 5 batches of particles, each executed in parallel over 5 processors, using 25 ranks in total. Spatial subcommunicators are used for domain decomposition algorithm for the iterative solvers involved in solving the forward equations for each particle, and ensemble subcommunicators are used to transfer particle states (and noises) during resampling.}
	\label{fig:enter-label}
\end{figure}

Our present implementation using Firedrake allows use to 
combine spatial domain decomposition for each particle
with ensemble parallelism across particles. The algorithms discussed above involve independent calculations for each particle, with the exception of the resampling step, when particle states (and noise increments) need to be replaced with copies from others.
Ensemble parallelism means dividing the ensemble of particles into batches and executing the independent calculations for each batch, before updating the particles from copies, which may come from other batches.
In our implementation, we use distributed memory parallelism using the Message Passing Interface (MPI) protocol, see Figure \ref{fig:enter-label}. 

Our present algorithm for resampling is quite naive: we just compute which particles need to be replaced by copies of which other ones, and send and receive from batches as necessary. Since communication of entire models states between ranks is costly (and likely to dominate the algorithm cost for large ensemble sizes), a more sophisticated approach should optimise the order of the particles after resampling to minimise communication.
Alternatively, algorithms such as the Islands Particle
Filter should be considered \cite{verge2015parallel}.
In the present work we do not investigate parallel performance. We just note that this combined parallelism
is possible in our code framework, and will present
a thorough investigation in future work on more challenging problems in 2D and 3D.

\section{Numerical investigations}\label{sec4}

In our numerical experiments of particle filtering applied to the SCH model, we use an ensemble of 150 particles. In the standard context, the velocity of the SCH model (\ref{fde1}) is observed every 5 time steps.
For all numerical experiments, we choose the length $L$ of the spatial domain to be 40. There are $N=100$ equispaced cells in the decomposition of the interval $[0,L]$, and the model time step is $\Delta t=0.025$. We explore multiple scenarios to evaluate stochastic filtering. In the initial scenario, observation data is gathered from the entire spatial domain, while in the second scenario, only half of the spatial domain is observed. All the numerical simulations are conducted using our general purpose particle filter library \cite{nudging} which is built upon
Firedrake \cite{FiredrakeUserManual}.

\subsection{Experiment 1: Full domain  observed}

In this experiment, we took measurements at $M = 81$ equispaced grid points in the interval $[0,L]$. 
The observations were perturbed with iid $\mathcal{N}(0,0.5)$ measurement errors.

\subsubsection{Initialization of particles and truth}
In this experiment the initialization of particles and the true solution of the model problem (\ref{fde1}) is constructed in the following way.

For particle $n=1,2,\ldots,N_p$, we solve the (finite element discretisation of the) following elliptic problem on the periodic domain $[0,L]$ with zero boundary,
\begin{equation}\label{init_solve}
\left\{
\begin{array}{ll}
      (I -  \nabla^{2}) U^{0,1}_{n}  = |W_{n}|,\quad n=1,2,\ldots,N_p,\\[2pt]
       (I -  \nabla^{2}) U^{0,j+1}_{n} = U^{0,j}_{n},\quad n=1,2,\ldots,N_p,\quad j=1,2,
\end{array}
\right.
\end{equation}
where $W_n$ is the DG0 function where each basis coefficient (i.e. constant cell value) is sampled
from $\N(0, h)$. We then calculate initial conditions for the particles as $u_{n}^{p,0} = ((\alpha_{n})^{2}U^{0,3}_{n}+(\beta_{n})^{2})$ for $ n=1,2,\ldots, N_p$, where 
random  parameters  $\alpha_{n},\beta_{n} \sim \N(0,1)$ for $n=1,2,\ldots,N_p$.
Then, the initial condition for the truth $u_{0}$ is sampled from the same distributions as the particles.

For the validation of various filtering procedures, we use the ensemble mean $l_2$-norm relative error (EMRE), the relative bias (RB) and the relative forecast ensemble spread (RES). These are defined as follows,
\[
EMRE(u^a , u^p):=  \dfrac{1}{N_p}\sum_{n=1}^{N_p}\dfrac{\|u^a - u_{n}^{p}\|_{2}}{\|u^a \|_{2}},
\]
\[
RB(u^a , u^p):= \dfrac{\|u^a - \overline{u}^{p}\|_{2}}{\|u^a \|_{2}}, \quad \mbox{ensemble mean}\,\, \overline{u}^{p} := \dfrac{1}{N_p}\sum_{n=1}^{N_p}u_{n}^{p},
\]
\[
RES(u^p, \overline{u}^{p}):= \dfrac{1}{N_p -1}\sum_{n=1}^{N_p}\dfrac{\|u_{n}^{p} - \overline{u}^{p}\|_{2}}{\|u^{a} \|_{2}}.
\]
These quantities are computed for each ensemble step.
The purpose of these statistics is demonstrate that
the particle filters are stable (or not).

\begin{figure}
    \centering
    \includegraphics[width=0.6\linewidth]{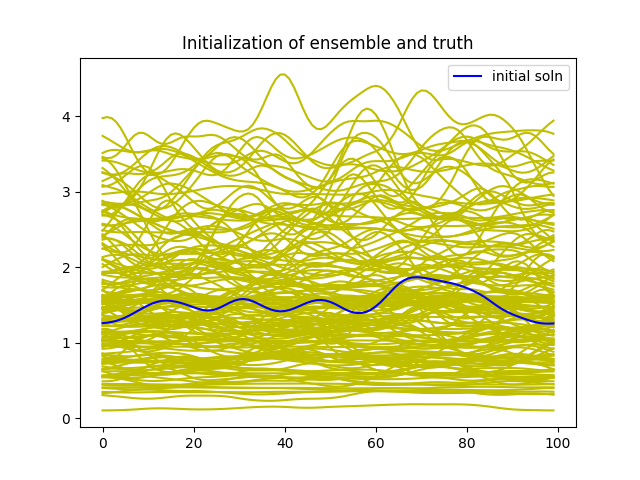}
    \caption{\small Initialization of all 150 particles and true state (the latter show in blue).}
    \label{Initiliazation}
\end{figure}

In Figure \ref{Initiliazation}, we have displayed the initialization of the particles and the true value.
We have presented the ensemble of particles in yellow colour in all the figures depicting the trajectories of the ensemble and true values. We observe that the cloud of particles is diverse.  The initial ensemble give good description of the initial uncertainty with $EMRE(u^a , u^{p,0}) = 0.48 0.12127201187663193
$ and $RB(u^a , u^{p,0}) = 0.12$.

\subsubsection{Bootstrap filter}
\begin{figure}
\subfloat[DA step 1]{\includegraphics[width=0.5\linewidth]{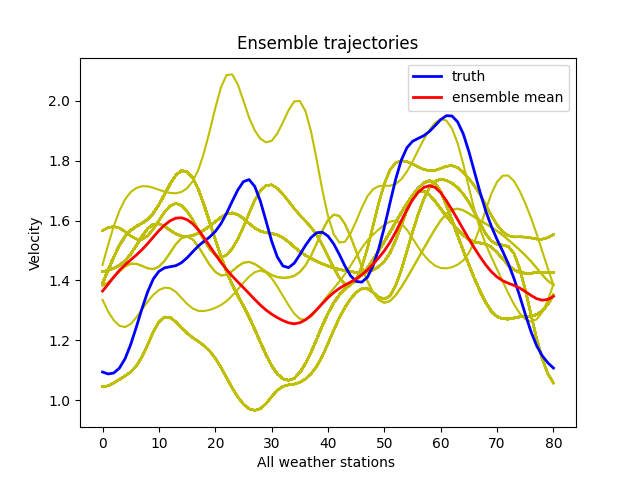}}
\hspace{-7mm}
\subfloat[DA step 100]{\includegraphics[width=0.5\linewidth]{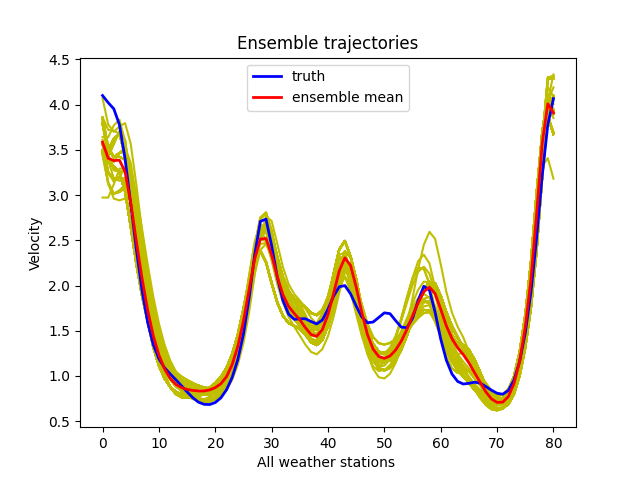}}\\
\subfloat[DA step 500]{\includegraphics[width=0.5\linewidth]{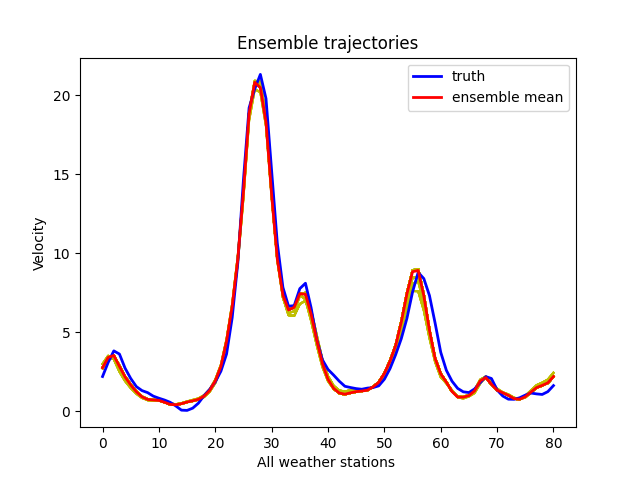}}\hspace{-7mm}
\subfloat[DA step 1000]{\includegraphics[width=0.5\linewidth]{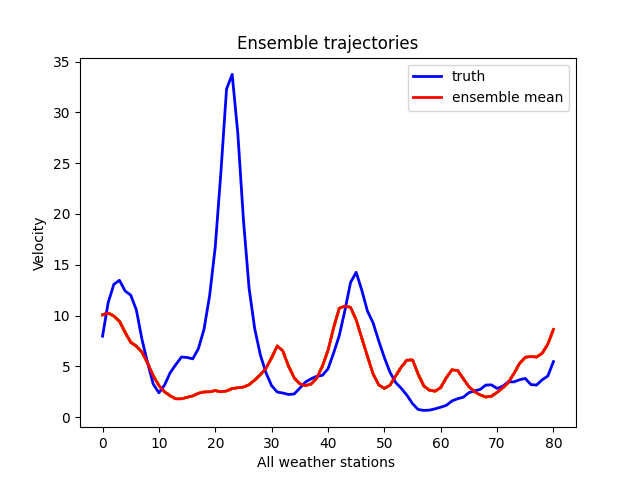}}\\
\caption{\small Comparison of the evolution of the true state  vs posterior ensemble and ensemble mean at data grids (weather stations). In order to assimilate
data we use bootstrap partile filter and outcome is displayed for the mentioned  assimilation steps.}
\label{SALT_bs_filter}
\end{figure}

\begin{figure}
\centering
\includegraphics[width=0.5\linewidth]{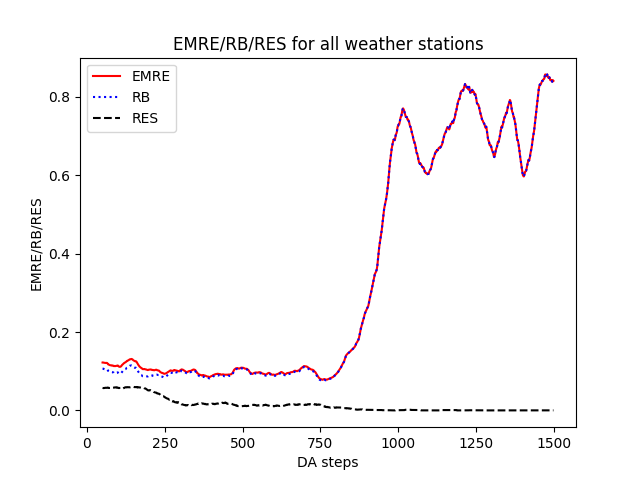}
\caption{\small Evolution of the relative ensemble mean error (EMRE), relative bias (RB) and ensemble spread (RES) associated with bootstrap particle filter.}
\label{SALT_bs_bias}
\end{figure}

Firstly, we discuss how the bootstrap filter performs for the SCH model. With the above initialization,  trajectories of the truth and particles are displayed for the different observation (data assimilation) steps in Figure \ref{SALT_bs_filter}. 
As time evolves, the ensemble spread reduces gradually but the cloud of particles does not track the truth. This phenomenon can be confirmed by Figure \ref{SALT_bs_filter}. The statistics of the difference between the truth and ensemble are discussed with the lenses of EMRE, RB and RES. These terms are calculated and displayed in Fig \ref{SALT_bs_bias} against assimilation steps. One can see the RES gradually decreases but the mean square error EMRE and relative bias RB diverge and saturate, indicating filter divergence.

\subsubsection{Tempering and jittering}
\begin{figure}
\subfloat[DA step 1]{\includegraphics[width=0.5\linewidth]{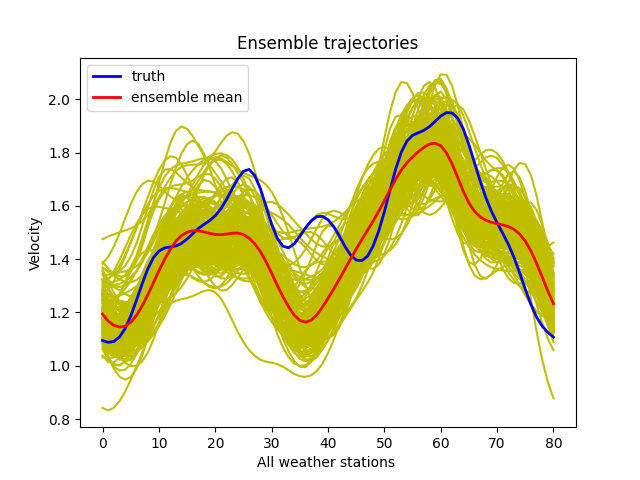}}
\hspace{-7mm}
\subfloat[DA step 100]{\includegraphics[width=0.5\linewidth]{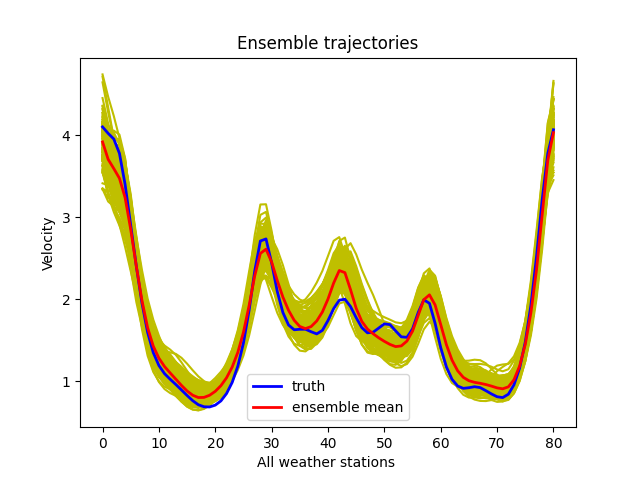}}\\
\subfloat[DA step 500]{\includegraphics[width=0.5\linewidth]{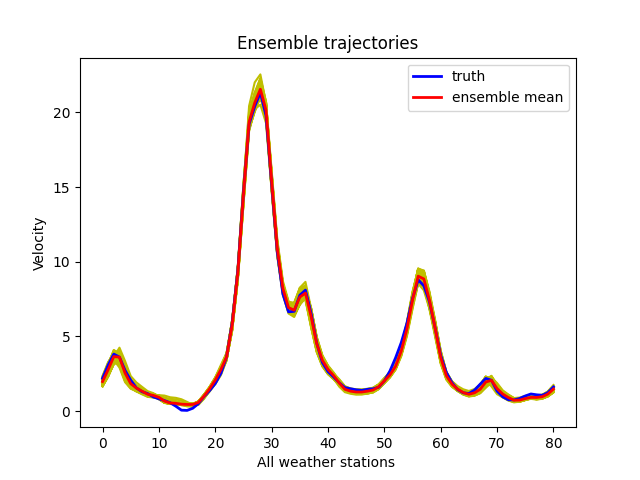}}\hspace{-7mm}
\subfloat[DA step 1000]{\includegraphics[width=0.5\linewidth]{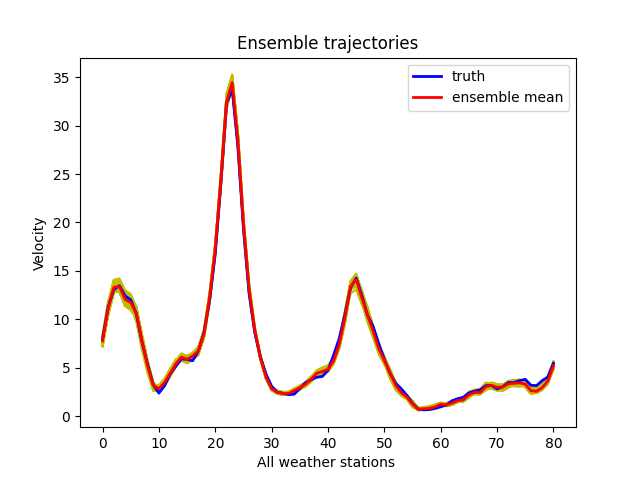}}\\
\caption{\small  Comparison of the evolution of the true state vs posterior ensemble and ensemble mean at data grids (weather stations). To assimilate
data we use {\it tempering} and {\it jittering} and the outcome is displayed for the mentioned assimilation steps. }
\label{SALT_temnp_filter}
\end{figure}

Now, we discuss the data assimilation algorithm that uses tempering and jittering. 
The jittering parameter $\delta$ is equal to $0.15$, and there are five jittering steps per tempering step.  Given that our framework allows for adapted tempering, the number of tempering steps needed for our numerical experiment falls within the range of 7 to 10.
In Fig. \ref{SALT_temnp_filter},  we exhibit a few instances to emphasize
the reduction of uncertainty resulting from applying the particle filter with tempering and jittering, as described in the previous section.
Additionally, it is possible to compare the ensemble after a single DA step using the tempering-jittering filter with the outcome from the bootstrap filter.
In particular, 
for the bootstrap filter, almost all
particles are replaced by duplicates of a small number
of particles, and the ensemble is not very diverse. However, the tempering allows the ensemble to be rediversified from those particles using the noise.

\begin{figure}
\centering
\includegraphics[width=0.5\linewidth]{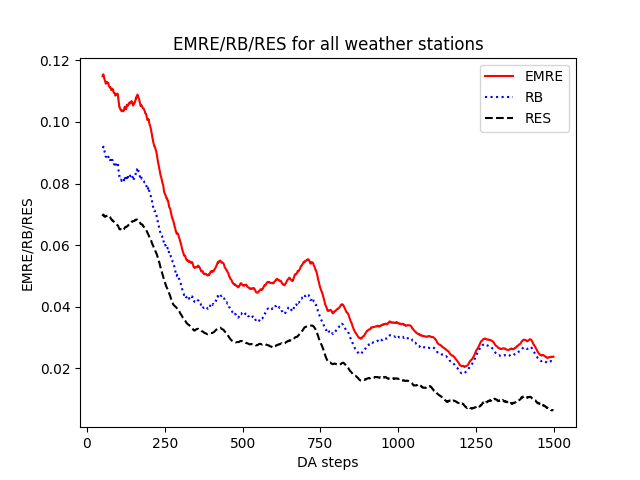}
\caption{\small Evolution of the relative ensemble mean error (EMRE), relative bias (RB) and ensemble spread (RES) associated with the filter using  {\it tempering} and {\it jittering} }
\label{SALT_temp_bias}
\end{figure}

We plot the EMRE, RB and RES in figure \ref{SALT_temp_bias}. 
One can see that
the EMRE and the RES are comparable and stable as time evolves. In contrast to the previous case, the particle filter follows the truth much better.

% {\bfseries: Maneesh: please can you make some comments about the number of tempering steps required.}

\subsubsection{Nudging}
\begin{figure}
\subfloat[DA step 1]{\includegraphics[width=0.5\linewidth]{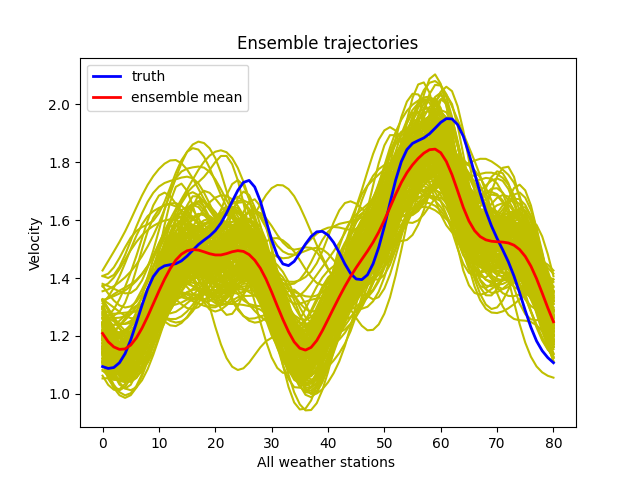}}
\hspace{-7mm}
\subfloat[DA step 100]{\includegraphics[width=0.5\linewidth]{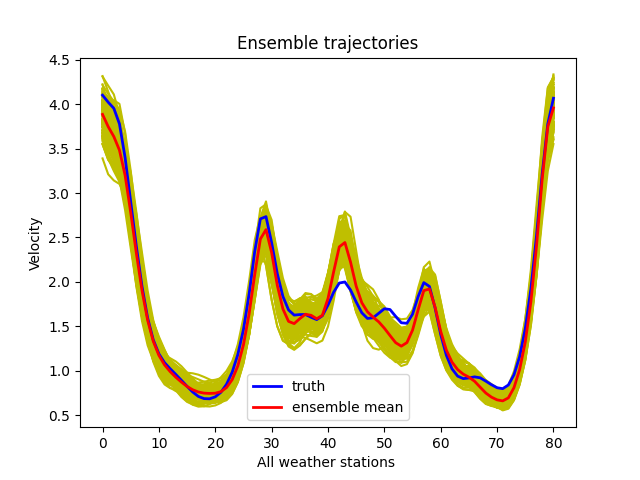}}\\
\subfloat[DA step 500]{\includegraphics[width=0.5\linewidth]{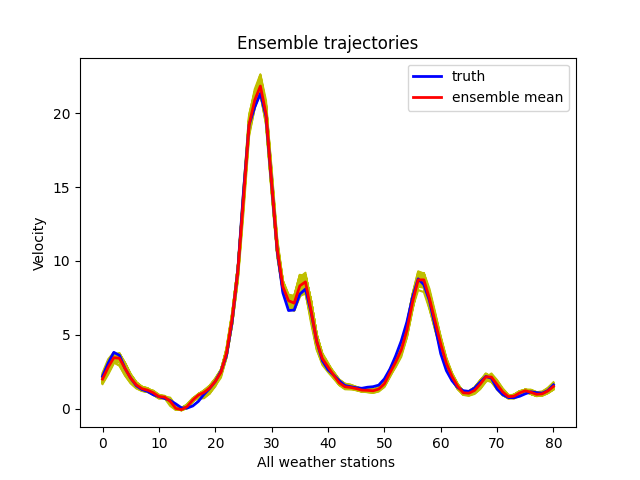}}\hspace{-7mm}
\subfloat[DA step 1000]{\includegraphics[width=0.5\linewidth]{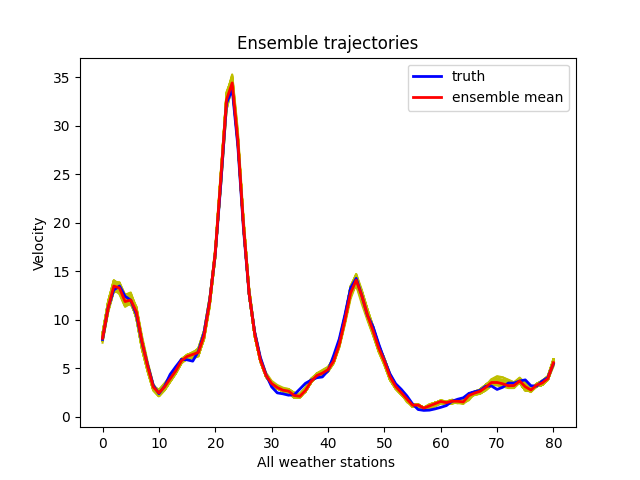}}\\
\caption{\small  Comparison of the evolution of the true state vs posterior ensemble and ensemble mean at data grids (weather stations). To assimilate
data we use  {\it nudging},  {\it tempering} and {\it jittering}  and outcome is displayed for the mentioned  assimilation steps }
\label{SALT_nudge_filter}
\end{figure}

\begin{figure}
\centering
\includegraphics[width=0.5\linewidth]{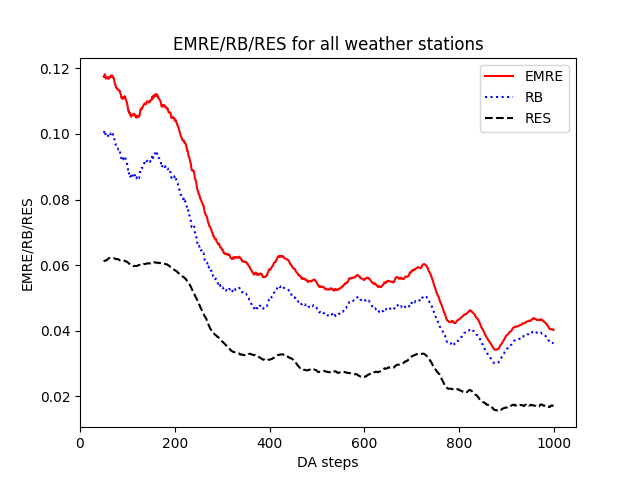}
\caption{\small Evolution of the relative ensemble mean error (EMRE), relative bias (RB) and ensemble spread (RES) associated with the filter using {\it nudging},  {\it tempering} and {\it jittering}}
\label{SALT_nudge_bias}
\end{figure}

We will now look at the performance of the data assimilation methodology which includes
nudging before using tempering and jittering.
In Figure \ref{SALT_nudge_filter},  we exhibit a few instances to emphasize
the reduction of uncertainty resulting from applying the particle filter with nudging, tempering and jittering. 
We plot the EMRE, RB and RES in Figure \ref{SALT_nudge_bias}. In these preliminary results
we observe that this particle filter is stable but does
not yet provide a dramatic improvement over the filter
without nudging.

\subsection{Experiment 2: Half domain  observed}
\begin{figure}
\subfloat[DA step 100]{\includegraphics[width=0.5\linewidth]{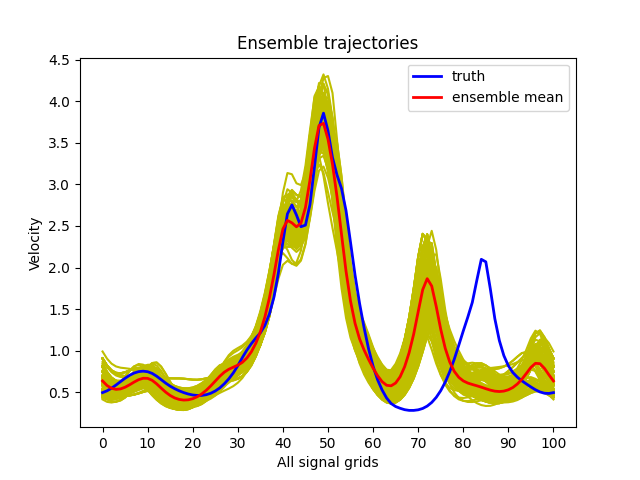}}
\hspace{-7mm}
\subfloat[DA step 500]{\includegraphics[width=0.5\linewidth]{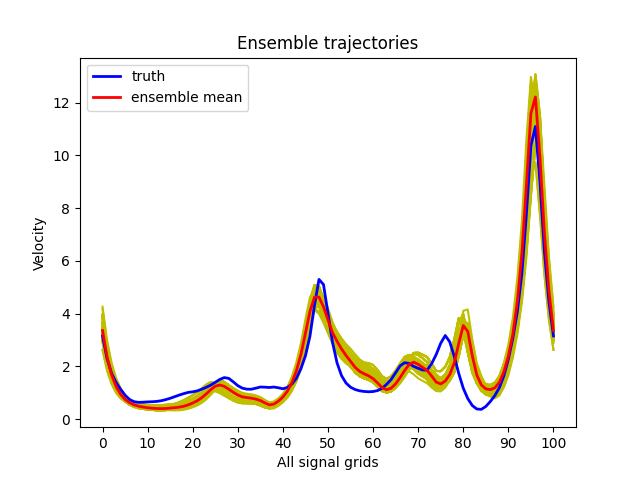}}\\
\subfloat[DA step 1000]{\includegraphics[width=0.5\linewidth]{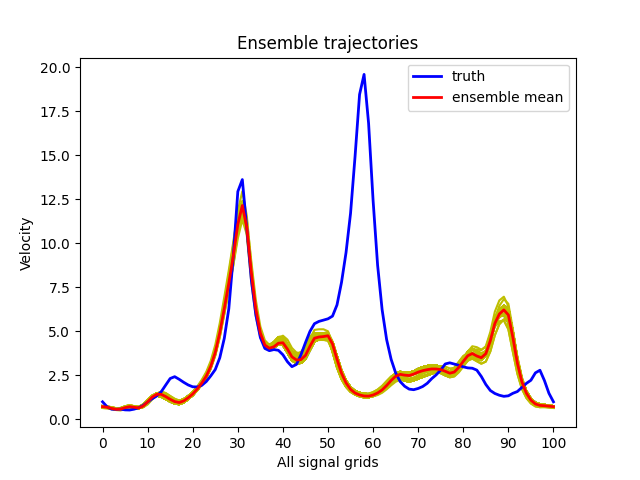}}\hspace{-7mm}
\subfloat[DA step 2000]{\includegraphics[width=0.5\linewidth]{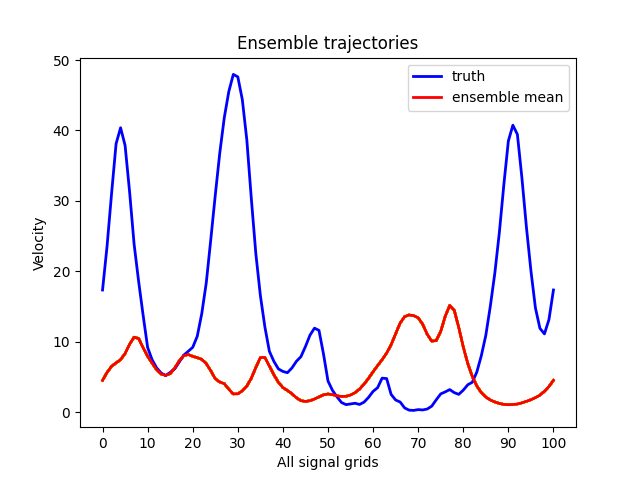}}\\
\caption{\small  Comparison of the evolution of the true state vs posterior ensemble and ensemble mean at data grids (weather stations). To assimilate
data we use bootstrap particle filter and outcome is displayed for the mentioned  assimilation steps }
\label{halfdomainSALT_bs}
\end{figure}

\begin{figure}
\subfloat[bootstrap filter]{\includegraphics[scale=.5]{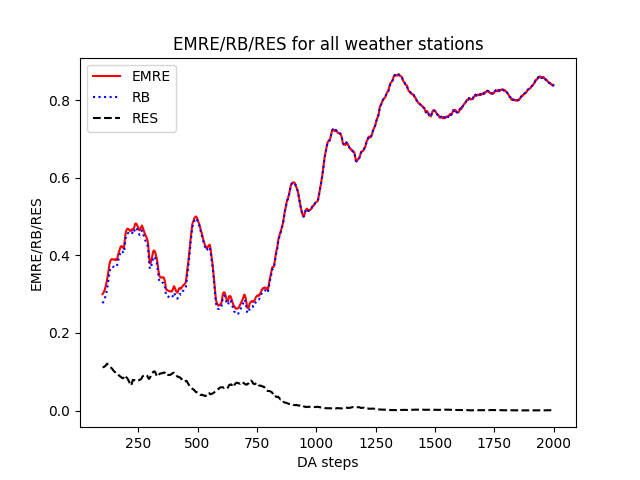}} \hspace{-5mm}
\subfloat[{\it tempering} and {\it jittering} ]{\includegraphics[scale=.5]{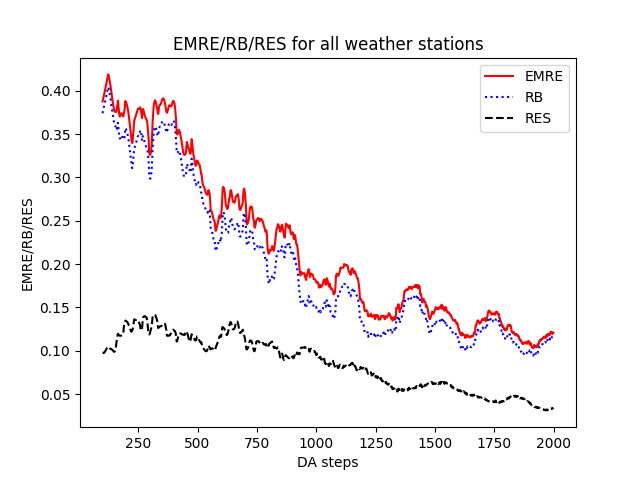}}
\caption{\small  Evolution of the relative ensemble mean error (EMRE), relative bias (RB) and ensemble spread (RES) for all signal grid points}
\label{halfdomainSALT_temp_bs_bias}
\end{figure}

In this experiment, the observation data is taken from half of the spatial domain, \emph{i.e}  $[0, L/2]$,
at $41$ equispaced points.
% Once more, we have chosen to use the same number of weather stations, \emph{i.e.}, $M=80$ as taken in the previous example. Specifically, 40 weather stations are within  $[0, L/2]$, and the remaining stations are located outside this domain.
With this modification, we now examine the bootstrap filter for the SCH model. From Figure \ref{halfdomainSALT_bs}, we see that truth is within the spread of the particles after initial data assimilation steps but particles lost track of truth even in the observed domain $[0, L/2]$. Also, we have plotted the values of EMRE, RB and RES in Figure \ref{halfdomainSALT_temp_bs_bias}, where one can observe that the error and bias (EMRE and RB)  significantly increase as time evolves.

\begin{figure}
\subfloat[DA step 100]{\includegraphics[width=0.5\linewidth]{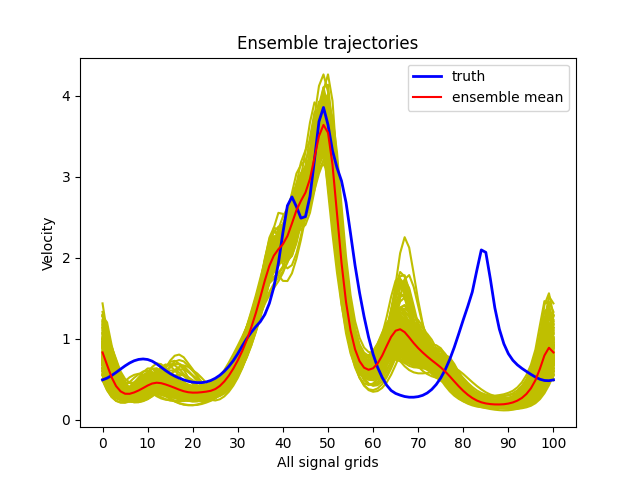}}
\hspace{-7mm}
\subfloat[DA step 500]{\includegraphics[width=0.5\linewidth]{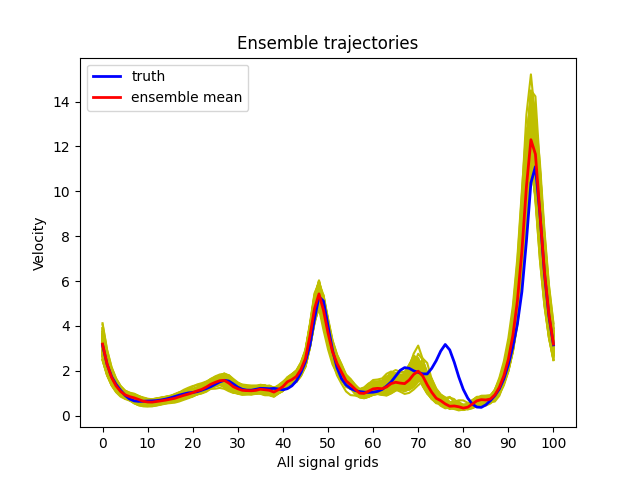}}\\
\subfloat[DA step 1000]{\includegraphics[width=0.5\linewidth]{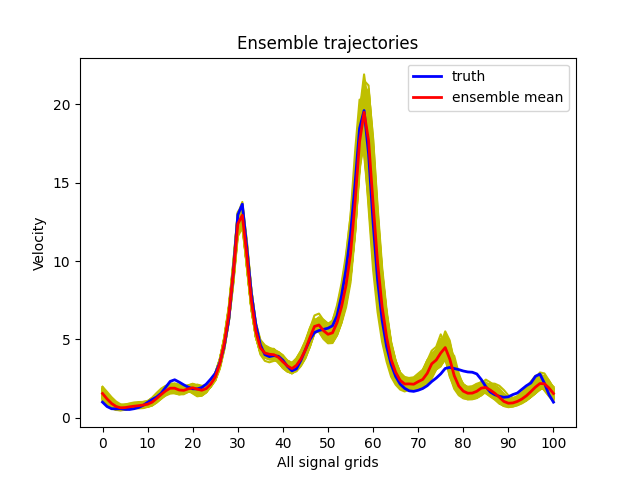}}\hspace{-7mm}
\subfloat[DA step 2000]{\includegraphics[width=0.5\linewidth]{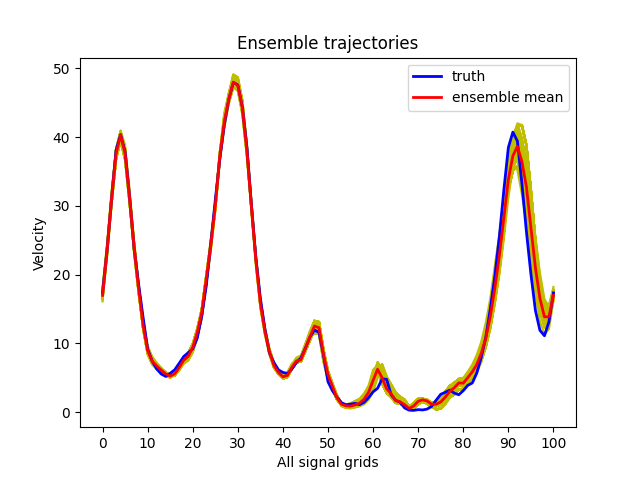}}\\
\caption{ Comparison of the evolution of the true state vs posterior ensemble and ensemble mean at data grids (weather stations). To assimilate
data we use  {\it tempering} and {\it jittering}  and outcome is displayed for the mentioned  assimilation steps }
\label{halfdomainSALT_temnp_filter}
\end{figure}

Next, we discuss the particle filter using tempering and jittering. We have used the same tempering-jittering procedure as discussed in the previous example. From Figure \ref{halfdomainSALT_temnp_filter}, we see that the ensemble surrounds the true value after some data assimilation steps in $[0, L/2]$. This does not occur 
initially; this is because it is not possible to reach all possible states of the posterior through modifications of the SALT noise. The ensemble does end up surrounding the true solution in the displayed results from the 100, 500, 1000 and 2000th data assimilation steps in this figure. This shows that
the particle filter is stable, and this is confirmed by
the EMRE, RB and RES measures. We leave the incorporation
of nudging for this example to further work.
% {\color{blue} we have old nudging results for half domain as well,  till 1000 DA steps, but no improvements compare to temper-jitter one}

\section{Concluding remarks}\label{sec5}
In this work, we investigated adaptive tempering, jittering and nudging techniques applied to the 
stochastic Camassa-Holm equation with SALT noise and
viscosity, demonstrating our new capability that
can be applied to arbitrary stochastic models written in Firedrake, made possible using MPI parallelism across
particles. The nudging filter involves solving nonlinear optimisation problems that are enabled with Firedrake's automated adjoint system. We demonstrated that these approaches lead to particle filters that are stable with relatively few (150) particles in cases where the classical bootstrap filter fails. 

In forthcoming work, we will undertake a detailed investigation of the accuracy of these filters when applied to 2D and 3D problems, making use of high performance computing. It is necessary to go beyond metrics such as ESS to properly determine the optimal value for $\delta$ and the optimal number of jittering steps; the gold standard is to compare against MCMC 
estimates of statistics.
We will also investigate whether the required accuracy can be more efficiently reached using Metropolis Adjusted Langevin (MALA) or Hybrid Monte Carlo (HMC) samplers in the jittering steps.

We will investigate the parallel performance of our implementation, and if necessary will develop more sophisticated ordering or subgrouping algorithms for resampling to achieve better parallel scalability.

Once we have established this capability, it is our goal to use data assimilation to calibrate the SALT parametrisation, modifying the Gaussian Mat\'ern field to model the differences between fine grid and coarse grid simulations, to facilitate faster data assimilation approaches.

\bibliographystyle{plain} 
\bibliography{sch_da} 
 
\end{document}